\documentclass[openany]{article}
\usepackage[margin=1in]{geometry} 
\usepackage{amsmath}

\renewcommand{\b}{\bf}
\newcommand{\norm}[1]{\left\lVert#1\right\rVert}

\newcommand{\cR}{\mathcal{R}}
\newcommand{\cG}{\mathcal{G}}
\newcommand{\cS}{\mathcal{S}}
\newenvironment{nproof}{\noindent\textit{Proof.}}{\hfill Q.E.D.}
\newcommand{\qed}{\hspace*{\fill}Q.E.D.} 

\newtheorem{theorem}{Theorem}[section]
\newtheorem{cor}[theorem]{Corollary}
\newtheorem{prop}[theorem]{Proposition}

\newtheorem{rem}[theorem]{Remark}

\newtheorem{lem}[theorem]{Lemma}
\newtheorem{example}[theorem]{Example}

\DeclareMathOperator{\Ric}{Ric}
\DeclareMathOperator{\wRic}{Ric_{a,b}}
\DeclareMathOperator{\PRic}{PRic}
\DeclareMathOperator{\vol}{Vol}
\DeclareMathOperator{\Hess}{Hess}
\DeclareMathOperator{\tr}{tr}

\title{On A Class of Weakly Weighted Einstein Metrics}
\author{ Zhongmin Shen\footnote{supported in part by a NSFC grant (12071423)} \ and \   Runzhong Zhao}
\date{}
\begin{document}
\maketitle
 
\begin{abstract} The notion of general weighted Ricci curvatures appears naturally in many problems. The N-Ricci curvature and the projective Ricci curvature are just two special  ones  with totally different geometric meanings. In this paper, we study general weighted Ricci curvatures. We find that Randers metrics of certain isotropic weighted Ricci curvature must have isotropic S-curvature. Then we classify  them  via their navigation expressions. We also find equations that characterize Randers metrics of almost isotropic weighted Ricci curvature. 
\end{abstract}

\section{Introduction}

	The $S$-curvature was first introduced  in the study of volume comparison in Finsler geometry (\cite{Sh1}), where a Bishop-Gromov type comparison theorem was proved. Since then  the $S$-curvature has been recognized 
	as an important non-Riemannian quantity in Finsler geometry.  Various combinations of Ricci curvature and $S$-curvature, now known as weighted Ricci curvatures, then came into 
	play in Finsler geometry. 
	
	In 2009, S. Ohta \cite{Oh} gave a diameter bound for a Finsler manifold $(M, F)$ with a volume form $dV$
under a positive lower bound of the $N$-Ricci curvature
	\[
		\Ric_N = \Ric + \dot{\mathbf{S}} - \frac1{N-n}\mathbf{S}^2,
	\]
where $\Ric$ is the Ricci curvature of $F$ and ${\bf S}$ is the S-curvature  of $(F, dV)$. 
However, in the extreme case when  $N = \infty$, a positive lower bound of 
	\[
		\Ric_\infty = \Ric + \dot{\mathbf{S}}
	\]
	does not imply the compactness of a Finsler manifold, but merely gives an upper bound for the volume.  This is first pointed out by S. Yin in \cite{Yin2}. Recently, X. Cheng and the first author give an estimate on the volume (\cite{ChSh4}).
On the other hand, the local expansion of the volume of small metric balls leads to the following weighted Ricci curvature \cite{Sh2}:
	\[
		\Ric + 3\dot{\mathbf{S}} - 3\mathbf{S}^2 .
	\]
There is another  quantity of this type occuring in the recent study of projective geometry, namely, the projective Ricci curvature 
	\[
		\PRic = \Ric + (n-1)\left[\frac{\dot{\mathbf{S}}}{n+1} + \frac{\mathbf{S}^2}{(n+1)^2}\right]
	\]
	was found to be a projective invariant \cite{Sh3}\cite{ShSu}. 

We now consider an important example. 
Let $h=\sqrt{h_{ij}(x)y^iy^j}$ be a Riemannian metric on an $n$-manifold $M$ satisfying
\begin{equation}
\Ric^h + a\Hess_h f - b\; df\otimes df = (n-1)\mu h^2,\label{Richf}
\end{equation}
where $a,b$ are constants and $\mu=\mu(x)$ is a scalar function on $M$. Let 
  $W= W^i\frac{\partial}{\partial x^i}$ be a vector field on $M$  satisfying
\begin{equation}
  W_{i;j}+ W_{j;i} =-4 c h_{ij},\label{Wh}
\end{equation}
where $W_i:=h_{ij}W^j$ and  $c=c(x)$ is a scalar function on $M$.  Assume that $\|W\|_h < 1$. We can define  a Randers metric  $F=\alpha+\beta$ by the {\it navigation data} $(h, W)$ as follows
\begin{equation}  
			F = \frac{\sqrt{\lambda h^2 + W_0^2}}{\lambda} -\frac{W_0}{\lambda },\label{navigation}
		\end{equation}
		where $\lambda := 1-\|W\|^2_h$. 
Let $ dV= e^{-f} dV_{BH}$  and ${\bf S}$ be the S-curvature of $(F, dV)$, where $dV_{BH}$ denote the Busemann-Hausdorff volume form. We shall show that  the Ricci curvature of $F$ satisfies
\begin{equation}
 \Ric + a\dot{\mathbf{S}} - b\mathbf{S}^2 = (n-1)\left(\frac{3\theta}{F} + \sigma\right)F^2 \label{RicS}
\end{equation}
where $\sigma =\sigma(x)$ is a scalar function and $\theta $ is a $1$-form on $M$. 
Both $\sigma$ and $\theta$ are expressed in terms of $f, c, \mu$,  $W$ and their derivatives with respect to $h$. Note that  (\ref{Wh}) is equivalent to that  ${\bf S}_{BH}=(n+1)c F$ and $\dot{\bf S}_{BH}= (n+1) c_0 F$. Thus if $f=0$, then  (\ref{RicS}) can be rewritten as
\begin{equation}
\Ric = (n-1) \Big \{  \frac{3\tilde{\theta}}{F} + \tilde{\sigma} \Big \}, \label{RicSt}
\end{equation}
where $ \tilde{\theta }:= \theta -\frac{(n+1)a}{3(n-1)} c_0   $ and $\tilde{\sigma}:=\sigma + \frac{(n+1)^2 bc^2}{n-1}  $.   This case has been discussed in \cite{ChSh3}.  However, when $f\not=0$, (\ref{RicS}) cannot be expressed as (\ref{RicSt}) in general. See Example \ref{ex1.3} below.

This leads to the notion of  an $(a, b)$-weighted Ricci curvature
\[ \wRic := \Ric + a\dot{\mathbf{S}} - b\mathbf{S}^2\]
where $a, b$ are constants. The $(a, b)$-weighted Ricci curvature can also be written in the following form
\begin{equation}
 \wRic = {\rm PRic} -\frac{\kappa}{n+1} \Big ( \dot{\b S}+\frac{4}{n+1} {\bf S}^2 \Big ) 
+\frac{\nu}{(n+1)^2} {\bf S}^2, \label{Rickv}
\end{equation}
where $ \kappa := (n-1)-a (n+1)$ and $\nu: = 3(n-1)-4a (n+1) -b(n+1)^2$. As we shall see,  there is a reason why we express $(a,b)$-weighted Ricci curvature in the form (\ref{Rickv}).

	A Finsler metric $F$ on an $n$-manifold $M$ is said to be \textit{weakly $(a, b)$-weighted Einstein} with respect to a volume form $dV=e^{-f}dV_{BH}$ if the $(a, b)$-weighted Ricci curvature satisfies
	\begin{equation}
		\wRic = (n-1)\left(\frac{3\theta}{F} + \sigma\right)F^2 \label{Einstein}
\end{equation}
	where 
	$\sigma$ is a scalar function and $\theta = \theta_iy^i$ is a one-form on $M$.  $F$ is said to be $(a,b)$-weighted Einstein with respect to a volume form $dV=e^{-f}dV_{BH}$ if it satisfies (\ref{Einstein}) with $\theta =0$. 

By the above definition, a Riemannian metric $h$ satisfying (\ref{Richf}) is $(a, b)$-weighted Einstein with respect  to $dV=e^{-f} dV_{h}$. 

To understand weakly $(a,b)$-weighted Einstein metrics, it is natural to study Randers metrics. 
 We first discover that when 
$\nu \ne 0$, the weakly $(a,b)$-weighted Einstein condition implies that 
the $S$-curvature is isotropic with respect to the Busemann-Hausdorff volume form.
	As a consequence, weakly $(a, b)$-weighted Einstein Randers metrics with $\nu\not=0$  can be obtained by navigations of homothetic field on weighted Einstein Riemannian metrics. 
	This result generalizes the results on weakly Einstein metrics in \cite{BaRo}, \cite{ChSh1} and \cite{ChSh3}. 

	\begin{theorem}\label{MTM}
		Let $F =\alpha+\beta$ be a Randers metric on an $n$-dimensional manifold $M$ 
		defined by a navigation $W$ on a Riemannian metric $h$
		Assume that $\nu\not=0$.
		$F$ is weakly $(a, b)$-weighted Einstein satisfying 
		\[
			\wRic = (n-1)\left(\frac{3\theta}{F} + \sigma\right)F^2
		\]
		with respect to a volume form $dV=e^{-f}dV_{BH}$ if and only if $h$ is  $(a,b)$-weighted Einstein satisfying (\ref{Richf}) with respect to $dV=e^{-f}dV_h$  and $W$ satisfies (\ref{Wh})
		for some scalar functions $f,c,\mu$ on $M$. 
		In this case 
		\begin{equation}\label{Meqs}
			\sigma = \mu - c^2 - 2c_iW^i + \frac1{n-1}\Big \{ -af_{i;j}W^iW^j + af_i\cS^i - bc^2(n+1)^2 + bf_if_jW^iW^j \Big \}
		\end{equation}
		and 
		\begin{equation}\label{Meqr}
			\theta_i = \frac1{3(n-1)}\Big \{\Big [3(n-1) + a(n+1)\Big ]c_i + 2af_{i;j}W^j + 2af_j\cS^j_{\;i} - 2cf_i\Big [ a + (n+1)b\Big ] - 2bf_if_jW^j \Big \} 
		\end{equation}
		where  $\cS^i$ and $\cS^j_{\;i} = h^{jk}\cS_{ki}$ are given in (\ref{eq14}) below. 
	\end{theorem}
	\begin{rem}{\rm 
		As shown below in proposition \ref{prop41}, the 'if' part in theorem \ref{MTM} does not rely on the assumption $\nu\ne 0$.
		A particularly interesting class of examples may arise from the case where $b = 0$, where the Ricci almost gradient solitons lead to Randers metrics with Einstein condition for $\Ric_\infty$. It is well-known that
		positive lower bound of $\Ric_\infty$ implies finite volume for forward-complete Finsler manifolds, which inspires the search for non-compact Randers metrics with positive lower bound of $\Ric_\infty$ constructed in this way. 
		}
	\end{rem}
	
	The following example shows that a weakly $(a,b)$-weighted Einstein Randers metric is not necessarily weakly Einstein. 
	\begin{example}\label{ex1.3}\rm
		Let $h = \sqrt{dt^2 + u(t)^2(dx^2 + dy^2)}$ on $\mathbf{R}^3$. Using O'Neil's formula (see also \cite{Be}) it's easy to see that $\Ric^h = -\frac{2u''(t)}{u(t)}dt^2 - (u(t)u''(t) + u'(t)^2)(dx^2
		+ dy^2)$. In case that $f = f(t)$ is a function of $t$ only, equation (\ref{Richf}) is equivalent to
		\begin{equation}\begin{aligned}\label{WEin}
			u(t)^2(af''(t) - bf'(t)^2) = u(t)u''(t) - u'(t)^2 + au(t)u'(t)f'(t)
		\end{aligned}\end{equation}
		where $\mu = -\frac{2u''(t)}{u(t)} + af''(t) - bf'(t)^2$. 

		An easy example comes when $u(t) = t$. Note that in particular $\Ric^h = -(dx^2 + dy^2)$ in this case, hence $h$ is not Einstein. Nevertheless one checks that 
		\[
		f(t)= \int_t^1\frac{1+\sqrt{2}\tanh(\sqrt{2}\ln s)}{s}ds
		\] satisfies equation (\ref{WEin}) with $a = b = 1$. This function $f$ is well-defined on the half space $\{t > 0\}$.

		On the other hand, there is an obvious killing field $W = -y\partial x + x\partial y$ whose norm $\norm{W}_h< 1$ in the region $\{t^2(x^2+y^2) < 1\}$. Thus in the
		region $\{(t,x,y)\mid t > 0, x^2 + y^2 < t^{-2}\}$, the navigation data $(h = \sqrt{dt^2 + t^2(dx^2 + dy^2)}, W = -y\partial x + x\partial y)$ defines a Randers metric which is weakly
		$(1,1)$-weighted Einstein.

		In view of Theorem 7.4.2 in \cite{ChSh1}, this Randers metric is not weakly Einstein since $h$ is not Einstein.
	\end{example}
	
	As a special case of our theorem, we obtain the following
	\begin{cor}
		Let $F=\alpha+\beta$ be a Randers metric on an $n$-dimensional manifold $M$ 
		defined by a navigation $W$ on a Riemannian metric $h$, 
		$F$ is weakly weighted Einstein satisfying 
		\[
			\Ric_\infty = (n-1)\left(\frac{3\theta}{F} + \sigma\right)F^2
		\]
		with respect to some volume form $dV$ if and only if $h$ is a Ricci almost gradient soliton $\Ric^h + \Hess_hf  = (n-1)\mu h^2$
and $\cR_{ij} = -2ch_{ij}$ for some
		scalar functions $c$. In this case we have $dV = e^{-f}dV_{BH}$, 
		\begin{equation}
			\sigma = \mu - c^2 - 2c_iW^i + \frac1{n-1}\left( f_i\cS^i -f_{i;j}W^iW^j \right)
		\end{equation}
		and 
		\begin{equation}
			\theta_i = \frac2{3(n-1)}\left[(2n - 1)c_i + f_{i;j}W^j + f_j\cS^j_{\;i} - cf_i\right]
		\end{equation}
	\end{cor}

One can refer to  \cite{Pi} for Ricci almost gradient solitons in Riemannian geometry.

	Complete characterization of weakly $(a,b)$-weighted Einstein Randers metric is more subtle when $\nu= 0$.  The $(a,b)$-weight Ricci curvature becomes
		\[\wRic = {\rm PRic} -\frac{\kappa}{n+1} \Big ( \dot{\b S}+\frac{4}{n+1} {\bf S}^2 \Big ) .\]
	In this case, the weakly $(a, b)$-weighted Einstein condition does not imply that the S-curvature is isotropic. Nevertheless, we can still obtain some equations on $\alpha$ and $\beta$ that characterize Randers metrics $F=\alpha+\beta$ of weakly $(a,b)$-weighted Ricci curvature. See Proposition \ref{prop5.1} below. If in addition $\kappa =0$, then the $(a, b)$-weighted Ricci curvature $\wRic = {\rm PRic}$ is just the projective Ricci curvature. We obtain a result that generalizes the result in \cite{ShSu} (see Proposition \ref{PRic} below).

\section{Preliminaries}
	Let $F = \alpha + \beta$ be a Randers metric on an $n$-manifold $M$, where $\alpha = \sqrt{a_{ij}(x)y^iy^j}$ is a Riemannian metric and $\beta = b_i(x)y^i$ is a 1-form with
	$\norm{\beta}_\alpha = \sqrt{a^{ij}(x)b_i(x)b_j(x)} < 1$. With a little abuse of notation, we shall denote by $;$ the covariant derivative with respect to the Riemannian metrics, so for $\alpha$ we write
	\[
		b_{i;j} = \dfrac{\partial b_i}{\partial x^j} - b_k{}^\alpha\Gamma^{k}_{ij}
	\]
	where ${}^\alpha\Gamma$ is the Christoffel symbol of $\alpha$. We shall also adopt the notations
	\begin{equation}\begin{aligned}
		r_{ij} : = \frac12\left(b_{i;j} + b_{j;i}\right), &\quad s_{ij} : = \frac12\left(b_{i;j} - b_{j;i}\right) \\
		s_j: = b^is_{ij},\quad r_j := b^ir_{ij},&\quad e_{ij} = r_{ij} + b_is_j + b_js_i
	\end{aligned}\end{equation}
	and further
	\begin{equation}\begin{aligned}
		q_{ij} := r_{im}s^m_{\;j},&\quad t_{ij} := s_{im}s^m_{\;j}\\
		q_j := b^iq_{ij} = r_ms^m_{\;j},&\quad t_j:= b^it_{ij} = s_ms^m_{\;j}
	\end{aligned}\end{equation}
	It is a well-known fact that the spray coefficients $G^i$ of $F$ and those $G^i_\alpha$ of $\alpha$ are related by 
	\begin{equation}
		G^i = G^i_\alpha + Py^i + Q^i
	\end{equation}
	where 
	\begin{equation}\label{Spray}
		P: = \frac{e_{00}}{2F} - s_0 = \frac{e_{ij}y^iy^j}{2F} - s_iy^i, \quad Q^i := \alpha s^i_{\;0} = \alpha s^i_{\;j}y^j
	\end{equation}
	with the convention that an index $0$ means the contraction with $y^i$.

	The spray $G:= y^i\frac{\partial}{\partial x^i} - 2G^i \frac{\partial}{\partial y^i}$  is a vector field on the tangent bundle, that 
	uniquely defines a family of transformation $\mathbf{R}_y := R^i_{\;k}dx^k\otimes\frac{\partial}{\partial x^i}: T_xM\to T_xM$. $\mathbf{R}$  is called the
	{\it Riemann curvature} of the spray $G$. In local coordinates, the coefficients $R^i_{\;k}$ are given by 
	\begin{equation}
		R^i_{\;k} = 2\frac{\partial G^i}{\partial x^k} - y^j\frac{\partial^2G^i}{\partial x^j\partial y^k} + 2G^j\frac{\partial^2G^i}{\partial y^j\partial y^k} - \frac{\partial
		G^i}{\partial y^j}\frac{\partial G^j}{\partial y^k}
	\end{equation}
	The trace of the Riemann curvature does not depend on the choice of local coordinates, and is called the {\it  Ricci curvature}:
	\begin{equation}
		\Ric := \tr(\mathbf{R}_y) = R^i_{\;i}
	\end{equation}
	For Randers metrics, the Ricci curvature is given by \cite{BaRo}\cite{ShYi}
	\begin{equation}
		\Ric = \Ric^\alpha + (2\alpha s^m_{\;0;m} - 2t_{00} - \alpha^2t^m_{\;m}) + (n-1)\Xi
	\end{equation}
	where $\Ric^\alpha$ is the Ricci curvature of $\alpha$ and
	\[
		\Xi = \frac{2\alpha}{F}(q_{00} - \alpha t_0) + \frac{3}{4F^2}(r_{00} - 2\alpha s_0)^2 - \frac1{2F}(r_{00;0} - 2\alpha s_{0;0})
	\]

	Alternatively, we will describe Randers metrics as solutions of the navigation problem on a Riemannian manifold. Let $h$
	a Riemannian metric and $W$ be a vector field with $\norm{W}_h < 1$. We may define a Finsler metric $F$ by 
	\begin{equation}
		\norm{\frac{y}{F(x,y)} - W}_h = 1
	\end{equation}
	It is well known that the Finsler metric $F$ obtained this way is a Randers metric. Indeed, $F = \alpha + \beta$
	with
	\begin{equation}\begin{aligned}
		a_{ij} =& \frac{h_{ij}}{\lambda} + \frac{W_iW_j}{\lambda^2}\\
		b_i =& -\frac{W_i}{\lambda}
	\end{aligned}\end{equation}
	where $\lambda = 1 - \norm{W}_h^2$.
	We shall denote
	\begin{equation}\begin{aligned} \label{eq14}
		\cR_{ij}:=\frac12\left(W_{i;j} + W_{j;i}\right),&\quad \cS_{ij} := \frac12\left(W_{i;j} - W_{j;i}\right)\\
		\cS_j:= W^i\cS_{ij},&\quad \cR_j := W^i\cR_{ij}, \quad \cR:= \cR_jW^j
	\end{aligned}\end{equation}
	where $;$ denotes the covariant derivative with respect to $h$, along
	with $\xi: = y - F(x,y)W$. By construction we have $h(x,\xi) = F(x,y)$.  

	Volume form comes into play when we introduce the weighted Ricci curvature using the S-curvature. In local coordinates, a volume form is described by a positive function $\sigma$:
	\begin{equation}
		dV = \sigma(x)dx^1\wedge\cdots\wedge dx^n
	\end{equation}
	The quantity 
	\[
		\tau(x,y) := \ln \frac{\sqrt{\det g_{ij}(x,y)}}{\sigma(x)}
	\]
	is called the \textit{distortion} and its rate of change along geodesics is measured by \textit{$S$-curvature}. Namely, let $c(t)$ be a geodesic with $c(0) = x$ and $\dot{c}(0) = y\in T_xM\setminus\{0\}$, we have 
	\begin{equation}\begin{aligned}
		\mathbf{S}(x,y) := \left.\frac{d}{dt}\right|_{t=0}\left[\tau(c(t), \dot{c}(t))\right],\quad\dot{\mathbf{S}}(x,y) := \left.\frac{d}{dt}\right|_{t=0}\left[\mathbf{S}(c(t), \dot{c}(t))\right]
	\end{aligned}\end{equation}
	In short, we have 
	\begin{equation}\begin{aligned}
		\mathbf{S} = \tau_{|i}y^i,\quad\dot{\mathbf{S}} = \mathbf{S}_{|i}y^i
	\end{aligned}\end{equation}
	where $|$ denotes the horizontal covariant derivative of $F$.

	Every Finsler metric $F$ induces a volume form called the {\it Busemann-Hausdorff volume form}. It is expressed in local coordinates by
	\[
		dV_{BH} :=\sigma_{BH}(x)dx^1\cdots dx^n =  \frac{\vol(\mathbf{B}^n(1))}{\vol{\{(y^i)\in\mathbf{R}^n\mid F(x,y) < 1\}}}dx^1\cdots dx^n
	\]
	
	A Finsler metric  $F$ is said of \textit{isotropic S-curvature} if $\mathbf{S} = (n+1)cF$ for some scalar function $c$.  The following lemma is well-known.
\begin{lem}\label{lemisoS}(\cite{ChSh2}\cite{Xi})
 For a Randers metric $F$ expressed by $\alpha+\beta$ or by the navigation data $(h, W)$, the following are equivalent for a scalar function $c=c(x)$, 

(a) $\mathbf{S} = (n+1)cF$ with respect to
	the Busemann-Hausdorff volume form ,

(b) $e_{00} = 2c(\alpha^2 - \beta^2)$,

(c) ${\cal R}_{00} = -2 c h^2$.

\end{lem}

	In the general case, let $dV = e^{(n+1)\phi}dV_\alpha$, one derives that \cite{ChSh1}
	\begin{equation}\begin{aligned}\label{S}
		\frac{\mathbf{S}}{n+1} = \frac{e_{00}}{2F} - (s_0 + \phi_0)
	\end{aligned}\end{equation}
	and using (\ref{Spray})
	\begin{equation}\begin{aligned}\label{Sdot}
		\frac{\dot{\mathbf{S}}}{n+1} =& \left[\frac{e_{00|l}}{2F}-(s_{0|l}+\phi_{0|l})\right]y^l\\
		=& \frac{e_{00;0}}{2F}-(s_{0;0}+\phi_{0;0}) - 2\left[\frac{e_{00}}{F} - (s_0+\phi_0)\right] \left(\frac{e_{00}}{2F} - s_0\right) - \frac{2\alpha}{F}e_{l0}s^l_{\;0} +2\alpha (s_l+\phi_l)s^l_{\;0}
	\end{aligned}\end{equation}

\section{Isotropic S-curvature}\label{section3}

Recall that the $(a,b)$-weighted Ricci curvature $\wRic=\Ric + a \dot{\bf S}- b {\bf S}^2$ can be expressed as
\[
 \wRic = {\rm PRic} -\frac{\kappa}{n+1} \Big ( \dot{\b S}+\frac{4}{n+1} {\bf S}^2 \Big ) 
+\frac{\nu}{(n+1)^2} {\bf S}^2,
\]
where $ \kappa := (n-1)-a (n+1)$ and $\nu: = 3(n-1)-4a (n+1) -b(n+1)^2$.
In this section, we shall study
 weakly $(a,b)$-weighted Einstein Randers metrics under the assumption $\nu\not=0$:
\[
			\wRic = (n-1)\left(\frac{3\theta}{F} + \sigma\right)F^2
		\]

	\begin{lem}\label{lem3.1}  Assume that $\nu\not=0$.
		If a  Randers metric on an $n$-manifold is 
		 weakly $(a, b)$-weighted Einstein, then it has isotropic $S$-curvature with respect to the Busemann-Hausdorff volume form.
	\end{lem}

	\begin{nproof}
		Taking $dV = e^{(n+1)\phi}dV_\alpha$, we have 
		\[\begin{aligned}
			0 =& F^2\wRic - (n-1)\left(\frac{3\theta}{F} + \sigma\right)F^4 \\=& F^2\Ric^\alpha + F^2(2\alpha s^m_{\;0;m} - 2t_{00} - \alpha^2t^m_{\;m}) + (n-1)F^2\Xi - (n-1)(3\theta F^3 + \sigma F^4)\\
			&+a(n+1)\bigg[\frac{F}{2}e_{00;0} - F^2(s_{0;0}+\phi_{0;0}) - e_{00}^2 + Fe_{00}(3s_0 + \phi_0) \\&- 2F^2s_0(s_0 + \phi_0) - 2\alpha Fe_{l0}s^l_{\;0} + 2\alpha F^2(s_l+\phi_l)s^l_{\;0}\bigg]-b(n+1)^2\left[\frac{e_{00}}{2} - F(s_0 + \phi_0)\right]^2.
		\end{aligned}\]
		The above equation can be  reorganized as follows:
		\begin{equation}\label{RicMain}\begin{aligned}
			0 =& \frac{1}{4}\nu e_{00}^2 +FP_1(y) + F^2P_2(y) + F^3P_3(y) + F^4P_4(y)
		\end{aligned}\end{equation}
		where $P_i$'s are polynomials with coefficients being functions of $x$ only. 		
Replacing $F$ by $\alpha+\beta$, we can rewrite   (\ref{RicMain}) as  
		\[
			A(y) + \alpha B(y) = 0,
		\]
where
		\begin{equation}
			A(y): = \frac{1}{4}\nu e_{00}^2 + \beta P_1 + (\alpha^2 + \beta^2)P_2 + (3\alpha^2\beta + \beta^3)P_3 + (\alpha^4 + 6\alpha^2\beta^2 + \beta^4)P_4  \label{A(y)}
\end{equation}
		and 
		\begin{equation}
\begin{aligned}
			B(y) := P_1 + 2\beta P_2 + (\alpha^2 + 3\beta^2)P_3 + (4\alpha^2\beta + 4\beta^3)P_4  \label{B(y)}.
		\end{aligned}
\end{equation}
Observe that $A$ and $B$ are polynomials in $y$, where $\alpha$ is irrational, hence $A(y) = B(y) = 0$. 
		
		Now that $A(y) - \beta B(y) = 0$, we have
		\begin{equation}\label{AB}
			\frac{1}{4}\nu  e_{00}^2 + (\alpha^2 - \beta^2)\left[P_2(y) + 2\beta P_3 + (\alpha^2 + 3\beta^2)P_4\right] = 0
\end{equation}
		Since $\norm{b}_\alpha < 1$, $\alpha^2 - \beta^2$ is an irreducible polynomial in $y$. Hence $\nu\ne 0$ leads to
		\begin{equation}\label{isotS}
			e_{00} = 2c(x)(\alpha^2 - \beta^2)
		\end{equation}
		for some scalar function $c=c(x)$, thus the  $S$-curvature is isotropic with respect to the Busemann-Hausdorff volume form. 
	\end{nproof}


\section{The case when $\nu\not=0$}

In this section, we are going to prove Theorem \ref{MTM}. 	
In Lemma \ref{lem3.1}, we have shown that if a Randers metric  $F=\alpha+\beta$  is weakly $(a, b)$-weighted Einstein with $\nu\not =0$, then (\ref{isotS}) holds or equivalently the S-curvature is isotropic with respect to the Busemann-Hausdorff volume form (see Lemma \ref{lemisoS} above).  Express $F=\alpha+\beta$  using a navigation data $(h, W)$ as in (\ref{navigation}).  Then the Riemann curvature of $F$ can be expressed  in terms of the Riemann curvature of $h$ and the covariant derivatives of $W$  (\cite{ChSh1}). Then we can derive equations in $h$ and $W$ that characterize weakly $(a, b)$-weighted Einstein  Randers metrics.

	\begin{prop}\label{prop41}
		Let $F$ be a Randers metric on an $n$-dimensional manifold $M$ defined by a navigation data $(h, W)$. Suppose that 
		the S-curvature with respect to Busemann-Hausdorff volume form $dV_{BH}$ is isotropic $\mathbf{S}_{BH} = (n+1)cF$. 
		Then $F$ is weakly $(a, b)$-weighted Einstein satisfying 
		\begin{equation}
			\wRic = (n-1)\left(\frac{3\theta}{F} + \sigma\right)F^2 \label{abRic}
		\end{equation}
		with respect to a volume form $dV= e^{-f} dV_{BH}$ 
		if and only if $\Ric^h + a\Hess_hf - b\; (df\otimes df) = (n-1)\mu h^2$ for some scalar functions $f$ and $\mu$. In this case,
		\begin{equation}
			\sigma = \mu - c^2 - 2c_iW^i + \frac1{n-1}\left[-af_{i;j}W^iW^j + af_i\cS^i - bc^2(n+1)^2 + bf_if_jW^iW^j \right]
		\end{equation}
		and 
		\begin{equation}
			\theta_i = \frac1{3(n-1)}\Big \{ \Big [ 3(n-1) + a(n+1)\Big ] c_i + 2af_{i;j}W^j + 2af_j\cS^j_{\;i} - 2cf_i\Big [a + (n+1)b\Big ] - 2bf_if_jW^j \Big \} 
		\end{equation}
	\end{prop}

	\begin{nproof}
		Let $G^i$ and $\cG^i$ be spray coefficients of $F$ and $h$, respectively. It is well-known that 
		\[\begin{aligned}
			\cG^i - G^i &= F\cS^i_{\;0} + \frac12F^2(\cR^i + \cS^i) - \frac12\left(\frac{y^i}{F} - W^i\right)(2F\cR_0 - \cR_{00} - F^2\cR).\end{aligned}\]
Recall that $F(x,y) = h(x, \xi) =: \tilde{h}$, where $\xi: = y - F(x,y)W$.
Since $\mathbf{S}_{BH} = (n+1)cF$, we have $\cR_{ij} = -2ch_{ij}$, and it follows that 
		\[
			2F\cR_0 - \cR_{00} - F^2\cR = -4c\tilde{h}W_i(\xi^i + \tilde{h}W^i) + 2ch_{ij}(\xi^i + \tilde{h}W^i)(\xi^j + \tilde{h}W^j) 
			+ 2c\tilde{h}^2\norm{W}_h^2 = 2c\tilde{h}^2
		\]
Then \[\cG^i - G^i= \tilde{h}\cS^i_{\;j}(\xi^j + \tilde{h}W^j) + \frac12\tilde{h}^2(\cR^i + \cS^i) - c\tilde{h}\xi^i
		\]
Now suppose that $F$ is weakly $(a, b)$-weighted Einstein as above, 
		taking a scalar function $f$ so that $dV = e^{-f}dV_{BH}$,  then $\dot{\mathbf{S}} = \dot{\mathbf{S}}_{BH} + \Hess_Ff = (n+1)c_0F 
		+ \Hess_Ff$. Then
		\[\begin{aligned}
			\Hess_Ff = f_{i|j}y^iy^j =& f_{i;j}(\xi^i + \tilde{h}W^i)(\xi^j + \tilde{h}W^j) + 2f_i(\cG^i - G^i)\\
			=& f_{i;j}(\xi^i + \tilde{h}W^i)(\xi^j + \tilde{h}W^j) + f_i\left(2\tilde{h}\cS^i_{\;j}(\xi^j + \tilde{h}W^j) + \tilde{h}^2(\cR^i + \cS^i)
			-2c\tilde{h}\xi^i\right)
		\end{aligned}\]
		By lemma 7.4.1 in \cite{ChSh1}, the weakly $(a,b)$-weighted Einstein condition (\ref{abRic})
		is equivalent to
		\begin{equation}\label{RicF}\begin{aligned}
			& \Ric^h(\xi) - (n-1)\tau\tilde{h}^2 + 3(n-1)c_i(\xi^i + \tilde{h}W^i)\tilde{h}\\
			=& (n-1)\left(3\theta_i(\xi^i + \tilde{h}W^i)\tilde{h} + \sigma\tilde{h}^2\right) - a(n+1)c_i(\xi^i + \tilde{h}W^i)\tilde{h}\\ &- a\left[f_{i;j}(\xi^i + \tilde{h}W^i)
			(\xi^j + \tilde{h}W^j) + f_i\left(2\tilde{h}\cS^i_{\;j}(\xi^j + \tilde{h}W^j) + \tilde{h}^2(\cR^i + \cS^i)-2c\tilde{h}\xi^i\right)\right]\\
			& + b\left[(n+1)c\tilde{h} + f_i(\xi^i + \tilde{h}W^i)\right]^2
		\end{aligned}\end{equation}
		where $\tau = c^2 + 2c_iW^i$. 
		
		Note that $\tilde{h}$ is irrational in $\xi$, separating rational and irrational terms in the above equation we have 
		\begin{equation}\label{1rat}\begin{aligned}
			\Ric^h - (n-1)\tau\tilde{h}^2 +& 3(n-1)c_iW^i\tilde{h}^2 = (n-1)(3\theta_iW^i + \sigma)\tilde{h}^2 - a(n+1)c_iW^i\tilde{h}^2 \\
			&-a\Hess_h{f}_{ij}\xi^i\xi^j - a\tilde{h}^2 f_{i;j}W^iW^j - 2af_i\tilde{h}^2\cS^i_{\;j}W^j - af^i\tilde{h}^2(\cR^i + \cS^i)\\
			&+bf_if_j\xi^i\xi^j + b\tilde{h}^2\left[(n+1)c + f_iW^i\right]^2
		\end{aligned}\end{equation}
		and
		\begin{equation}\label{1irrat}\begin{aligned}
			3(n-1)c_i\xi^i\tilde{h} =& 3(n-1)\theta_i\xi^i\tilde{h} - a(n+1)c_i\xi^i\tilde{h} \\
			-&2af_{i;j}\xi^iW^j\tilde{h} - 2af_j\cS^j_{\;i}\xi^i\tilde{h} + 2acf_i\xi^i\tilde{h} + 2bf_i\xi^i\left[(n+1)c
			+ f_jW^j\right]\tilde{h}
		\end{aligned}\end{equation}

		Now suppose that $F$ is weakly $(a,b)$-weighted Einstein. From (\ref{1rat}) we obtain
		\[
			\Ric^h(\xi) + a\Hess_hf(\xi) - bf_if_j\xi^i\xi^j = (n-1)\mu \tilde{h}^2
		\]
		where 
		\[\begin{aligned}
			\mu =& \tau + \sigma + 3\theta_iW^i \\&- \frac1{n-1}\left[((3+a)n + (a-3))c_iW^i + af_{i;j}W^iW^j + 2af_i\cS^i_{\;j}W^j + af_i(\cR^i + \cS^i)
			- b((n+1)c + f_iW^i)^2\right]
		\end{aligned}\]
		
		Conversely, suppose that $\Ric^h +a\Hess_hf - bdf\otimes df = (n-1)\mu h^2$	and $\cR_{ij} = -2ch_{ij}$.
		According to equation (\ref{1irrat}), we may choose 
		\begin{equation}\label{theta}
			\theta_i = \frac1{3(n-1)}\Big \{ (3(n-1) + a(n+1))c_i + 2af_{i;j}W^j + 2af_j\cS^j_{\;i} - 2cf_i[a + (n+1)b] - 2bf_if_jW^j \Big \} 
		\end{equation}
		and then
		\[
			\sigma = \mu - c^2 - 2c_iW^i + \frac1{n-1}\left[-af_{i;j}W^iW^j + af_i\cS^i - bc^2(n+1)^2 + bf_if_jW^iW^j \right]
		\]
		It is  easy to check that equation (\ref{1rat}) also holds and it follows that 
		\[
			\wRic = (n-1)\left(\frac{3\theta}{F} + \sigma\right)F^2
		\]
		with respect to $dV = e^{-f}dV_{BH}$. 
	\end{nproof}

\bigskip

\noindent{\it Proof of Theorem \ref{MTM} }: 
Theorem \ref{MTM} now follows by  combining Lemma \ref{lem3.1} with Proposition \ref{prop41}.  \qed

\bigskip

	We may also find equations on $\alpha$ and $\beta$  that characterize weakly $(a, b)$-weighted Einstein metrics. We shall continue to use the notations in Section \ref{section3}. 

Assume that $F=\alpha+\beta$ is weakly $(a, b)$-weighted Einstein with  $\nu\ne 0$.  By Lemma \ref{lem3.1},  we have 
\begin{equation}
e_{00} = 2c(x)(\alpha^2 - \beta^2),\label{e_00}
\end{equation}
Then the $P_i$'s in (\ref{RicMain}) are given by
	\[\begin{aligned}
		P_1 =& (\alpha^2 - \beta^2)\Big\{\kappa (4\beta c^2 + 2c(s_0 + 3\phi_0) - c_0) - 2c\nu(s_0 + \phi_0)\Big\}\\
		P_2 =& \Ric^\alpha - 2\beta s^m_{\;0;m} - 2t_{00} - \beta^2t^m_{\;m} - a(n+1)[2\beta \phi_ls^l_{\;0} + \phi_{0;0} - \phi_0^2] \\
			& + \kappa \left[2\beta t_0 + s_{0;0} - 4\beta cs_0 + s_0^2 - 3(s_0 + \phi_0)^2\right] + \nu(s_0 + \phi_0)^2\\
		P_3 =& 2\beta t^m_{\;m} + 2s^m_{\;0;m} + 2a(n+1)\phi_ls^l_{\;0} - 3(n-1)\theta - 2\kappa t_0\\
		P_4 =& -t^m_{\;m} - (n-1)\sigma
	\end{aligned}\]
	Now the equation  $B(y) - \frac{2\beta}{\alpha^2 - \beta^2}(A(y) - \beta B(y)) = 0$ gives 
	\[
		P_1(y) + (\alpha^2 - \beta^2)(P_3(y) + 2\beta P_4(y) - 2\beta c^2\nu) = 0
	\]
	from which we obtain
	\begin{equation}\label{eqs1}\begin{aligned}
		s^m_{\;0;m} =& (n-1)(3\theta/2 + \beta\sigma - \phi_ls^l_{\; 0}) \\
		&+  \kappa \Big \{\frac{c_0}{2} - 2\beta c^2 - c(s_0+3\phi_0) + t_0+\phi_ls^l_{\;0}\Big \} + c\nu\left(\beta c + s_0 + \phi_0\right)
	\end{aligned}\end{equation}
	Plugging into (\ref{AB}) and solving for $\Ric^\alpha$ we have
	\begin{equation}\label{eqr1}\begin{aligned}
		&\Ric^\alpha = 2t_{00} + \alpha^2 t^m_{\;m}  +(n-1)[ \phi_{0;0}-\phi_0^2+ 3\beta\theta +\sigma (\alpha^2 + \beta^2)]\\
		&- \kappa \Big \{ \phi_{0;0}-\phi_0^2 + s_0^2 - 3(s_0 + \phi_0)^2 + s_{0;0} + \beta c_0 - 4\beta^2c^2 - 6\beta c(s_0 + \phi_0)\Big \}\\
		&-\nu(s_0 + \phi_0 + \beta c)^2 - \nu c^2\alpha^2
	\end{aligned}\end{equation}

Conversely, if (\ref{e_00}),  (\ref{eqs1}) and 
		(\ref{eqr1}) hold, then it is easy to see that $F$ is weakly $(a, b)$-weighted Einstein.

We have proven the following.

	\begin{theorem} Let $a, b$ be two constants satisfying $\nu \ne 0$ and $F = \alpha + \beta$ be a Randers metric on an $n$-dimensional manifold $M$. 
		$F$ is weakly $(a, b)$-weighted Einstein satisfying 
		\[
			\wRic = (n-1)\left(\frac{3\theta}{F} + \sigma\right)F^2
		\]
		with respect to some volume form $dV=e^{(n+1)\phi}dV_{\alpha}$  if and only if equations (\ref{e_00}),  (\ref{eqs1}) and 
		(\ref{eqr1}) are satisfied with some scalar functions $\phi$ and $c$. 
	\end{theorem}

\section{ The case when $\nu = 0$ and $\kappa \not=0$ }

In this section we shall consider the $(a, b)$-weighted Ricci curvature when 
$\nu  = 0$. In this case, the $(a, b)$-Ricci curvature is given by 
\[ \wRic = {\rm PRic} -\kappa  \Big [ \frac{\dot{\b S}}{n+1} +\frac{4 {\bf S}^2 }{(n+1)^2}\Big ],\]
where $ \kappa :=  (n-1)-a(n+1)$. 
We are going to derive an equivalent condition for the weakly $(a, b)$-weighted Einstein condition
	\begin{equation}
		\wRic = (n-1)\left(\frac{3\theta}{F} + \sigma\right)F^2.\label{Ricnu=0}
	\end{equation}

 We shall continue to use the notations in Section \ref{section3}. 

Assume that $F$ is weakly $(a,b)$-weighted Einstein satisfying (\ref{Ricnu=0}). 
First note that equation (\ref{AB}) now gives 
	\begin{equation}\label{eq1}
		P_2 + 2\beta P_3 + (\alpha^2 + 3\beta^2)P_4 = 0
	\end{equation}
	Plugging back into either $A(y) = 0$ or $B(y) = 0$ we have 
	\begin{equation}\label{eq2}
		P_1 + (\alpha^2 - \beta^2)P_3 + 2\beta(\alpha^2 - \beta^2)P_4 = 0
	\end{equation}
In this case,  the polynomials $P_i$'s can be simplified a little:
	\[\begin{aligned}
		P_1 =& -\kappa\left[\frac12e_{00;0} + 2\beta e_{l0}s^{l}_{\;0} - e_{00}(s_0 + 3\phi_0)\right]\\
		P_2 =& \Ric^\alpha - 2\beta s^m_{\;0;m} - 2t_{00} - \beta^2t^m_{\;m} - (n-1)[\phi_{0;0} - \phi_0^2+2\beta \phi_ls^l_{\;0} ] \\
			& + \kappa\left[ \phi_{0;0} - \phi_0^2+2\beta \phi_ls^l_{\;0} +2e_{l0}s^{l}_{\;0} + s_0^2 + 2\beta t_0 - 3(s_0 + \phi_0)^2 + s_{0;0}\right]\\
		P_3 =& 2\beta t^m_{\;m} + 2s^m_{\;0;m} + 2a(n+1)\phi_ls^l_{\;0} - 3(n-1)\theta - 2\kappa t_0\\
		P_4 =& -t^m_{\;m} - (n-1)\sigma
	\end{aligned}\]
Solving for $s^m_{\;0;m}$ using equation (\ref{eq2}) we have 
	\begin{equation}\label{eqs}\begin{aligned}
		s^m_{\;0;m} =& (n-1)(3\theta/2 + \beta\sigma - \phi_ls^l_{\; 0}) \\
		&+  \kappa\Big \{ \dfrac{e_{00;0} + 4\beta e_{l0}s^{l}_{\;0} - 2e_{00}(s_0 + 3\phi_0)}{2(\alpha^2 - \beta^2)}+t_0 + \phi_ls^l_{\;0}\Big \}
	\end{aligned}\end{equation}
Plugging  into equation (\ref{eq1}) and solving for $\Ric^\alpha$,  we have 
	\begin{equation}\label{eqr}\begin{aligned}
		&\Ric^\alpha = 2t_{00} + \alpha^2 t^m_{\;m}  +(n-1)[ \phi_{0;0}-\phi_0^2+ 3\beta\theta +\sigma (\alpha^2 + \beta^2)]\\
		-& \kappa\Big \{ \phi_{0;0}-\phi_0^2 + s_0^2 - 3(s_0 + \phi_0)^2 + s_{0;0}+ 2 e_{l0}s^l_{\; 0}+\beta
\dfrac{e_{00;0} +4\beta e_{l0}s^l_{\; 0} - 2e_{00}(s_0 + 3\phi_0) }{\alpha^2 - \beta^2}\Big \}
	\end{aligned}\end{equation}

Now we assume that the dimension $n \geq 3$ and 
$\kappa \not=0$. 
It follows from (\ref{eqs}) that there is a $1$-form $\eta$ such that 
\begin{equation}
e_{00;0} + 4\beta e_{l0}s^{l}_{\;0} - 2e_{00}(s_0 + 3\phi_0) =2\eta (\alpha^2-\beta^2).\label{e000}
\end{equation}
In fact, (\ref{eqs}) allows us to write this $\eta$ as
\begin{equation}
\eta = \frac{1}{\kappa}\Big \{ s^m_{\;0;m} - (n-1)(3\theta/2 + \beta\sigma-\phi_ls^l_{\;0}\Big \} - t_0 - \phi_ls^l_0,\label{eqsm0m}
\end{equation}
and (\ref{eqr}) is reduced to
\begin{equation}\label{eqre}\begin{aligned}
		&\Ric^\alpha = 2t_{00} + \alpha^2 t^m_{\;m}  +(n-1)[ \phi_{0;0}-\phi_0^2+ 3\beta\theta +\sigma (\alpha^2 + \beta^2)]\\
		-& \kappa [\phi_{0;0}-\phi_0^2 +s_0^2 -3(s_0+\phi_0)^2 +s_{0;0} + 2 e_{0l}s^l_{\; 0} +2\beta \eta  ]
	\end{aligned}\end{equation}

Conversely, if $\alpha$ and $\beta$ satisfy (\ref{e000}) and 
		(\ref{eqre}) are satisfied with some scalar function $\phi$ and $1$-form $\eta$ defined by (\ref{eqsm0m}), then it is easy to verify that (\ref{Ricnu=0}) holds. 
We have proved the following

	\begin{prop} \label{prop5.1} Let $a, b$ be two constants satisfying $\nu  = 0$ and $\kappa \not=0$. 
		Let $F = \alpha + \beta$ be a Randers metric on an $n$-dimensional manifold $M$. 
		$F$ is weakly $(a, b)$-weighted Einstein satisfying 
		\[
			\wRic = (n-1)\left(\frac{3\theta}{F} + \sigma\right)F^2
		\]
		with respect to some volume form $dV$  if and only if equations (\ref{e000}) and 
		(\ref{eqre}) are satisfied with some scalar function $\phi$ and $1$-form $\eta$ defined by (\ref{eqsm0m}).
	\end{prop}

\begin{rem}{\rm 
	In case that $F$ is of isotropic $S$-curvature with respect to the Busemann-Hausdorff volume form, repeatedly using equation (\ref{isotS}) we obtain 
	\[
		\eta = c_0 - 4c^2\beta - 2c(s_0 + 3\phi_0), \quad e_{l0}s^{l}_{\;0} = -2\beta cs_0
	\]
	and the above condition reduces to (\ref{eqs1}) and (\ref{eqr1}) with $\nu = 0$.
		It is, however, unknown (to the authors) whether the existence of the 1-form $\eta$ (\ref{e000}) implies that the $S$-curvature is isotropic (with respect to the Busemann-Hausdorff volume form).
		}
	\end{rem}

\section{The projective Ricci curvature}

In the case when $\nu = 0$ and $\kappa=0$, 
the $(a, b)$-weighted Ricci curvature is just the projective Ricci curvature studied in \cite{ShSu}. 
\[ {\rm Ric}_{a, b}= {\rm PRic}.\]
In this case, equations (\ref{eqs}) and (\ref{eqr}) are reduced to
\begin{equation}
s^m_{\;0;m} = 
(n-1) (  3\theta/2 + \beta\sigma -\phi_ls^l_{\;0}),
 \label{eqs0}
\end{equation}
\begin{equation}
\Ric^\alpha =  2t_{00} + \alpha^2 t^m_{\;m} + (n-1)\left(\phi_{0;0} - \phi_0^2 + 3\beta\theta + \sigma (\alpha^2 + \beta^2)\right)
.\label{eqr0}
\end{equation}
Therefore  we obtain the following proposition due to Gasemnezhad-Rezaei-Gabrani.
	
	\begin{prop}\label{PRic} {(\cite{GRG})}
		Let $F = \alpha + \beta$ be a Randers metric on an $n$-dimensional manifold $M$,
		$F$ is weakly $(a, b)$-weighted Einstein satisfying 
		\begin{equation}
\PRic = (n-1)\left(\frac{3\theta}{F} + \sigma\right)F^2 \label{PRicWE}
\end{equation}
		with respect to some volume form $dV=e^{(n+1)\phi}dV_{\alpha}$, if and only if  equations (\ref{eqs0})and (\ref{eqr0}) are satisfied for some scalar function $\phi$.  
	\end{prop}
We shall remark that the equation (\ref{PRicWE}) is not a projective condition. Namely, if $F_1$ and $F_2$ are projectively related and $F_1$ satisfies (\ref{PRicWE}) for some volume form $dV$, then $F_2$ might not satisfy (\ref{PRicWE}) for any volume form. That is the reason why the authors in \cite{ShSu} study projectively Ricci-flat Finsler metrics. Anyway we 
 recovered the  result in \cite{ShSu} for projectively Ricci-flat Randers metric and Tabatabaeifar and collaborators' result\cite{Ta} for weighted projectively Ricci-flat metric by setting $\theta = 0$ and $\sigma = 0$. 
	
\bibliographystyle{plain}

\vskip 8mm

\noindent
Zhongmin Shen \& Runzhong Zhao \\
Department of Mathematical Sciences \\
Indiana University-Purdue University Indianapolis \\
IN 46202-3216, USA  \\
E-mail: zshen@math.iupui.edu\\
E-mail: runzzhao@iu.edu

\end{document}